\documentclass[11pt,a4paper,reqno]{amsart}  
\usepackage{amssymb} \usepackage{amscd} \usepackage{times}
\newcommand{\be} {\begin{equation}}
\newcommand{\ee} {\end{equation}}
\newcommand{\bea} {\begin{eqnarray}}
\newcommand{\eea} {\end{eqnarray}}
\newcommand{\Bea} {\begin{eqnarray*}}
\newcommand{\Eea} {\end{eqnarray*}}

\def\zbb{\mathbb{Z}}  
  
  \def\phi{\varphi}
 \def\p1{{\mathbb{P}^1_\zbb}}

\newtheorem{Theorem}{\quad Theorem}[section]

\begin{document}
\title{A compactness result for a system with weight and boundary singularity.}
\author{Samy Skander Bahoura} 
\address{Departement de Mathematiques, Universite Pierre et Marie Curie, 2 place Jussieu, 75005, Paris, France.}
\email{samybahoura@yahoo.fr, samybahoura@gmail.com}
\maketitle
\begin{abstract}
We give blow-up behavior for solutions to an elliptic system with Dirichlet condition, and, weight and boundary singularity. Also, we have a  compactness result for this elliptic system with regular H\"olderian weight and boundary singularity and Lipschitz condition. 
\end{abstract}


\section{Introduction and Main Results}

We set $ \Delta = \partial_{11} + \partial_{22} $  on open analytic domain $ \Omega $ of $ {\mathbb R}^2 $.

\bigskip

We consider the following equation:

$$ (P)   \left \{ \begin {split} 
       -\Delta u & = |x|^{2 \beta} V e^{v} \,\, &\text{in} \,\, & \Omega  \subset {\mathbb R}^2, \\
      - \Delta v & = W e^{u} \,\, &\text{in} \,\, & \Omega  \subset {\mathbb R}^2, \\
                    u & = 0  \,\,             & \text{in} \,\,    &\partial \Omega,\\ 
                    v & = 0  \,\,             & \text{in} \,\,    &\partial \Omega.                        
   \end {split}\right.
$$

Here, we assume that:

$$ 0 \leq V \leq b_1 < + \infty, \,\, e^u \in L^1({\Omega})\,\, {\rm and} \,\,  u \in W_0^{1,1}(\Omega), $$

$$ 0 \leq W \leq b_2 < + \infty, \,\, |x|^{2\beta} e^v \in L^1({\Omega})\,\, {\rm and} \,\,  v \in W_0^{1,1}(\Omega), $$

and,

$$  0 \in \partial \Omega, \,\,\, \beta \geq 0. $$
 
When $ u=v $ and $ \beta =0 $, the above system is reduced to an equation which was studied by many authors, with or without  the boundary condition, also for Riemann surfaces,  see [1-17],  one can find some existence and compactness results, also for a system.

Among other results, we  can see in [6] the following important Theorems ($ \beta=0 $):,

\smallskip

{\bf Theorem A.}{\it (Brezis-Merle [6])}.{\it  Consider the case of one equation; if $ (u_i)_i=(v_i)_i $ and $ (V_i)_i=(W_i)_i $ are two sequences of functions relatively to the problem $ (P) $ with, $ 0 < a \leq V_i \leq b < + \infty $, then, for all compact set $ K $ of $ \Omega $,

$$ \sup_K u_i \leq c = c(a, b, K, \Omega). $$}

\smallskip

{\bf Theorem B} {\it (Brezis-Merle [6])}.{\it Consider the case of one equation and assume that $ (u_i)_i $ and $ (V_i)_i $ are two sequences of functions relatively to the previous problem $ (P) $ with, $ 0 \leq V_i \leq b < + \infty $, and,

$$ \int_{\Omega} e^{u_i} dy  \leq C, $$

then, for all compact set $ K $ of $ \Omega $,

$$ \sup_K u_i \leq c = c(b, C, K, \Omega). $$}

Next, we call energy the following quantity:

$$ E= \int_{\Omega} e^{u_i} dy. $$

The boundedness of the energy is a necessary condition to work on the problem $ (P) $ as showed in $ [6] $, by the following counterexample ($\beta =0$):

\smallskip

{\bf Theorem C} {\it (Brezis-Merle [6])}.{\it  Consider the case of one equation, then there are two sequences $ (u_i)_i $ and $ (V_i)_i $ of the problem $ (P) $ with, $ 0 \leq V_i \leq b < + \infty $, and,

$$ \int_{\Omega} e^{u_i} dy  \leq C, $$

and

$$ \sup_{\Omega}  u_i \to + \infty. $$}

When $ \beta = 0 $, the above system have many properties in the constant and the Lipschitzian cases. Indeed we have (when $ \beta = 0 $):

In [12], Dupaigne-Farina-Sirakov proved (by an existence result of Montenegro, see [16]) that the solutions of the above system when $ V $ and $ W $ are constants can be extremal and this condition imply the boundedness of the energy and directly the compactness.
Note that in [11], if we assume (in particular) that $ \nabla \log V $ and $ \nabla \log W $ and $ V >a >0 $ or $ W > a'>0 $ and $ V, W $ are nonegative and uniformly bounded then the energy is bounded and we have a compactness result.

Note that in the case of one equation (and $ \beta = 0 $), we can prove by using the Pohozaev identity that if $ + \infty >b \geq V \geq a >0 $, $ \nabla V $ is uniformely Lipschitzian that the energy is bounded when $ \Omega $ is starshaped. In [15] Ma-Wei, using the moving-plane method showed that this fact is true for all domain $ \Omega $ with the same assumptions on $ V $. In [11] De Figueiredo-do O-Ruf extend this fact to a system by using the moving-plane method for a system.

Theorem C, shows that we have not a global compactness to the previous problem with one equation, perhaps we need more information on $ V $ to conclude  to the boundedness of the solutions. When $ \nabla \log V $ is  Lipschitz function and $ \beta = 0 $, Chen-Li and Ma-Wei see [7] and [15], showed that we have a compactness on all the open set. The proof is via the moving plane-Method of Serrin and Gidas-Ni-Nirenberg. Note that in [11], we have the same result for this system when $ \nabla \log V $ and $ \nabla \log W $ are uniformly bounded. We will see below that for a system we also have a compactness result when $ V $ and $ W $ are Lipschitzian and $ \beta \geq 0 $.

Now consider the case of one equation. In this case our equation have nice properties.

If we assume $ V $ with more regularity, we can have another type of estimates, a $ \sup + \inf $ type inequalities. It was proved by Shafrir see [17], that, if $ (u_i)_i, (V_i)_i $ are two sequences of functions solutions of the previous equation without assumption on the boundary and, $ 0 < a \leq V_i \leq b < + \infty $, then we have the following interior estimate:

$$ C\left (\dfrac{a}{b} \right ) \sup_K u_i + \inf_{\Omega} u_i \leq c=c(a, b, K, \Omega). $$

\bigskip

Now, if we suppose $ (V_i)_i $ uniformly Lipschitzian with $ A $ the
Lipschitz constant, then, $ C(a/b)=1 $ and $ c=c(a, b, A, K, \Omega)
$, see [5]. 
\smallskip

Here we are interested by the case of a system of this type of equation. First, we give the behavior of the blow-up points on the boundary, with weight and boundary singularity, and in the second time we have a proof of  compactness of the solutions to Gelfand-Liouville type system with weight and boundary singularity and Lipschitz condition.

\smallskip

Here, we write an extention of Brezis-Merle Problem (see [6]) to a system:

\smallskip

{\bf Problem}. Suppose that $ V_i \to  V $ and  
$ W_i \to  W $ in $ C^0( \bar \Omega ) $, with, $ 0 \leq V_i $ and $ 0 \leq W_i $. Also, we consider two sequences of solutions $ (u_i), (v_i) $ of $ (P) $ relatively to $ (V_i), (W_i) $ such that,

$$ \int_{\Omega} e^{u_i} dx \leq C_1,\,\,\,  \int_{\Omega} |x|^{2\beta} e^{v_i} dx \leq C_2,  $$

is it possible to have:

$$ ||u_i||_{L^{\infty}}\leq C_3=C_3(\beta, C_1, C_2, \Omega) ? $$

and,

$$ ||v_i||_{L^{\infty}}\leq C_4=C_4(\beta, C_1, C_2, \Omega) ? $$

In this paper we give a caracterization of the behavior of the blow-up points on the boundary and also a proof of the compactness theorem when $ V_i $ and $ W_i $ are uniformly Lipschitzian and $ \beta \geq 0 $. For the behavior of the blow-up points on the boundary, the following condition are enough,

$$ 0 \leq  V_i \leq b_1, \,\,\, 0 \leq  W_i \leq b_2, $$

The conditions $ V_i \to  V $ and $ W_i \to  W $ in $ C^0(\bar \Omega) $ are not necessary.

\bigskip

But for the proof of the compactness  for the system, we assume that:

$$ ||\nabla V_i||_{L^{\infty}}\leq  A_1, \,\,\, ||\nabla W_i||_{L^{\infty}}\leq  A_2,\,\, \beta \geq 0. $$

Our main result are:

\begin{Theorem}  Assume that $ \max_{\Omega} u_i \to +\infty $ and $ \max_{\Omega} v_i \to +\infty $ Where $ (u_i) $ and $ (v_i) $ are solutions of the probleme $ (P) $ with ($ \beta \geq 0 $), and:
 
 $$ 0 \leq V_i \leq b_1,\,\,\, {\rm and } \,\,\, \int_{\Omega}  e^{u_i} dx \leq C_1, \,\,\, \forall \,\, i, $$
 
and,

 $$ 0 \leq W_i \leq b_2,\,\,\, {\rm and } \,\,\, \int_{\Omega}  |x|^{2 \beta} e^{v_i} dx \leq C_2, \,\,\, \forall \,\, i, $$

 then;  after passing to a subsequence, there is a finction $ u $,  there is a number $ N \in {\mathbb N} $ and $ N  $ points $ x_1, x_2, \ldots, x_N \in  \partial \Omega $, such that, 

$$ \int_{\partial \Omega} \partial_{\nu} u_i  \phi  \to \int_{\partial \Omega} \partial_{\nu} u  \phi +\sum_{j=1}^N \alpha_j  \phi(x_j), \, \alpha_j \geq 4\pi, $$
$ \text{for\, any}\,\, \phi\in C^0(\partial \Omega) $, 
and,

$$ u_i \to u \,\,\, {\rm in }\,\,\, C^1_{loc}(\bar \Omega-\{x_1,\ldots, x_N \}). $$

$$ \int_{\partial \Omega} \partial_{\nu} v_i  \phi  \to \int_{\partial \Omega} \partial_{\nu} v  \phi +\sum_{j=1}^N \beta_j  \phi(x_j), \, \beta_j \geq 4\pi, $$
$ \text{for\, any}\,\, \phi\in C^0(\partial \Omega) $,  and,

$$ v_i \to v \,\,\, {\rm in }\,\,\, C^1_{loc}(\bar \Omega-\{x_1,\ldots, x_N \}). $$

\end{Theorem} 

 In the following theorem, we have a proof for the global a priori estimate which concern the problem $ (P) $.

\bigskip

\begin{Theorem}Assume that $ (u_i), (v_i) $ are solutions of $ (P) $ relatively to $ (V_i), (W_i) $ with the following conditions:

$$ x_1=0 \in \partial \Omega, \,\,\beta \geq 0, $$

and,

$$ 0 \leq  V_i \leq b_1, \,\,  ||\nabla V_i||_{L^{\infty}} \leq A_1,\,\, {\rm and } \,\,\, \int_{\Omega} e^{u_i} \leq C_1, $$

$$ 0 \leq  W_i \leq b_2, \,\,  ||\nabla W_i||_{L^{\infty}} \leq A_2,\,\, {\rm and } \,\,\, \int_{\Omega} |x|^{2\beta} e^{v_i} \leq C_2, $$

We have,

$$  || u_i||_{L^{\infty}} \leq C_3(b_1, b_2, \beta, A_1, A_2, C_1, C_2, \Omega), $$

and,

$$  || v_i||_{L^{\infty}} \leq C_4(b_1, b_2, \beta, A_1, A_2, C_1, C_2, \Omega), $$

\end{Theorem} 

\section{Proof of the theorems} 

\bigskip

\underbar {\it Proof of theorem 1.1:} 

\bigskip

We have:

$$ u_i, v_i \in W_0^{1,1}(\Omega). $$

Since $ e^{u_i} \in L^1(\Omega) $ by the corollary 1 of Brezis-Merle's paper (see [6]) we have $ e^{v_i} \in L^k(\Omega) $ for all $ k  >2 $ and the elliptic estimates of Agmon and the Sobolev embedding (see [1]) imply that:

$$ u_i \in W^{2, k}(\Omega)\cap C^{1, \epsilon}(\bar \Omega). $$ 

And,

We have:

$$ v_i, u_i \in W_0^{1,1}(\Omega). $$

Since $ |x|^{2\beta} e^{v_i} \in L^1(\Omega) $ by the corollary 1 of Brezis-Merle's paper (see [6]) we have $ e^{u_i} \in L^k(\Omega) $ for all $ k  >2 $ and the elliptic estimates of Agmon and the Sobolev embedding (see [1]) imply that:

$$ v_i \in W^{2, k}(\Omega)\cap C^{1, \epsilon}(\bar \Omega). $$

Since $ |x|^{2\beta} V_ie^{v_i} $ and $ W_ie^{u_i} $ are bounded in $ L^1(\Omega) $, we can extract from those two sequences two subsequences which converge to two nonegative measures $ \mu_1 $ and $ \mu_2 $. (This procedure is similar to the procedure of Brezis-Merle, we apply corollary 4 of Brezis-Merle paper, see [6]).

If $ \mu_1(x_0) < 4 \pi $, by a Brezis-Merle estimate for the first equation, we have $ e^{u_i} \in L^{1+\epsilon} $ around $ x_0 $, by the elliptic estimates, for the second equation, we have $ v_i \in W^{2, 1+\epsilon} \subset L^{\infty} $ around $ x_0 $, and , returning to the first equation, we have $ u_i \in L^{\infty} $ around $ x_0 $.

If $ \mu_2(x_0) < 4 \pi $, then $ u_i $ and $ v_i $ are also locally bounded around $ x_0 $.

Thus, we take a look to the case when, $ \mu_1(x_0) \geq 4 \pi $ and $ \mu_2(x_0) \geq 4 \pi $. By our hypothesis, those points $ x_0 $ are finite.

We will see that inside $ \Omega $ no such points exist. By contradiction, assume that, we have $ \mu_1(x_0) \geq 4 \pi $. Let us consider a ball $ B_R(x_0) $ which contain only $ x_0 $ as nonregular point. Thus, on $ \partial B_R(x_0) $, the two sequence $ u_i $ and $ v_i $ are uniformly bounded. Let us consider:

$$ \left \{ \begin {split} 
       -\Delta z_i & = |x|^{2\beta}V_i e^{v_i} \,\, &\text{in} \,\, & B_R(x_0)  \subset {\mathbb R}^2, \\
                       z_i & = 0  \,\,             & \text{in} \,\,    &\partial B_R(x_0).                                            
   \end {split}\right.
$$

By the maximum principle we have:

$$ z_i \leq u_i $$ 

and $ z_i \to z $ almost everywhere on this ball, and thus,

$$ \int e^{z_i} \leq \int e^{u_i} \leq C, $$

and,

$$ \int e^z \leq C.$$

but, $ z  $ is  a solution in $ W_0^{1,q}(B_R(x_0)) $, $ 1\leq q <2 $, of the following equation:

$$ \left \{ \begin {split} 
       -\Delta z & = \mu_1\,\, &\text{in} \,\, & B_R(x_0) \subset {\mathbb R}^2, \\
                     z & = 0  \,\,             & \text{in} \,\,    &\partial B_R(x_0).                                            
   \end {split}\right.
$$

with, $ \mu_1 \geq 4 \pi $ and thus, $ \mu_1 \geq 4\pi \delta_{x_0} $ and then, by the maximum principle in $ W_0^{1,q}(B_R(x_0)) $:

$$ z \geq -2 \log |x-x_0|+ C $$

thus,

$$ \int e^z = + \infty, $$

which is a contradiction. Thus, there is no nonregular points inside  $ \Omega $

Thus, we consider the case where we have nonregular points on the boundary, we use two estimates:

$$ \int_{\partial \Omega} \partial_{\nu} u_i d\sigma \leq C_1,\,\,\, \int_{\partial \Omega} \partial_{\nu} v_i d\sigma \leq C_2, $$
 
 and,
 
 $$ ||\nabla u_i||_{L^q} \leq C_q, \,\,\, ||\nabla v_i||_{L^q} \leq C'_q, \,\,\forall \,\, i\,\, {\rm and  }  \,\, 1< q < 2. $$

We have the same computations, as in the case of one equation.

We consider a points $ x_0 \in \partial \Omega $ such that:

$$ \mu_1(x_0) < 4 \pi. $$

We consider a test function on the boundary $ \eta $ we extend $ \eta $ by a harmonic function on $ \Omega $, we write the equation:

$$ -\Delta ((u_i-u)\eta) =|x|^{2\beta}(V_i e^{v_i}-Ve^v)\eta+ <\nabla (u_i-u)|\nabla \eta> = f_i $$

with,

$$ \int |f_i| \leq 4 \pi-\epsilon +o(1) < 4\pi-2\epsilon <4\pi, $$

$$ -\Delta ((v_i-v)\eta) =(W_i e^{u_i}-We^u)\eta+ <\nabla (v_i-v)|\nabla \eta> = g_i, $$

with,

$$ \int |g_i| \leq 4 \pi-\epsilon +o(1) < 4\pi-2\epsilon <4\pi, $$

By the Brezis-Merle estimate, we have uniformly, $ e^{u_i} \in L^{1+\epsilon} $ around $ x_0 $, by the elliptic estimates, for the second equation,  we have $ v_i \in W^{2, 1+\epsilon} \subset L^{\infty} $ around $ x_0 $, and , returning to the first equation, we have $ u_i \in L^{\infty} $ around $ x_0 $.

We have the same thing if we assume:

$$ \mu_2(x_0) < 4 \pi. $$

Thus, if $ \mu_1(x_0) < 4 \pi $ or $ \mu_2(x_0) < 4 \pi $, we have for $ R >0 $ small enough:

$$ (u_i,v_i) \in L^{\infty}(B_R(x_0)\cap \bar \Omega). $$

By our hypothesis the set of the points such that:

$$ \mu_1(x_0)  \geq  4 \pi,  \,\,\, \mu_2(x_0)   \geq 4 \pi, $$

is finite, and, outside this set $ u_i $ and $ v_i $ are locally uniformly bounded. By the elliptic estimates, we have the $ C^1 $ convergence to  $ u $ and $ v $  on each compact set of $ \bar \Omega- \{x_1, \ldots x_N\} $.

Indeed,

By the Stokes formula we have, 

$$ \int_{\partial \Omega} \partial_{\nu} u_i d\sigma \leq C, $$

We use the weak convergence in the space of Radon measures to have the existence of a nonnegative Radon measure $ \mu_1 $ such that,

$$ \int_{\partial \Omega} \partial_{\nu} u_i \phi  d\sigma \to \mu_1(\phi), \,\,\, \forall \,\,\, \phi \in C^0(\partial \Omega). $$

We take an $ x_0 \in \partial \Omega $ such that, $ \mu_1({x_0}) < 4\pi $. For $ \epsilon >0 $ small enough set $ I_{\epsilon}= B(x_0, \epsilon)\cap \partial \Omega $ on the unt disk or one can assume it as an interval. We choose a function $ \eta_{\epsilon} $ such that,

$$ \begin{cases}
    
\eta_{\epsilon} \equiv 1,\,\,\,  {\rm on } \,\,\,  I_{\epsilon}, \,\,\, 0 < \epsilon < \delta/2,\\

\eta_{\epsilon} \equiv 0,\,\,\, {\rm outside} \,\,\, I_{2\epsilon }, \\

0 \leq \eta_{\epsilon} \leq 1, \\

||\nabla \eta_{\epsilon}||_{L^{\infty}(I_{2\epsilon})} \leq \dfrac{C_0(\Omega, x_0)}{\epsilon}.

\end{cases} $$

We take a $\tilde \eta_{\epsilon} $ such that,

\begin{displaymath}  \left \{ \begin {split} 
      -\Delta \tilde \eta_{\epsilon}  & = 0              \,\, &&\text{in} \!\!&&\Omega \subset {\mathbb R}^2, \\
                  \tilde\eta_{\epsilon} & =  \eta_{\epsilon}   \,\,             && \text{in} \!\!&&\partial \Omega.               
\end {split}\right.
\end{displaymath}

{\bf Remark:} We use the following steps in the construction of $ \tilde \eta_{\epsilon} $:

We take a cutoff function $ \eta_{0} $ in $ B(0, 2) $ or $ B(x_0, 2) $:

1- We set $ \eta_{\epsilon}(x)= \eta_0(|x-x_0|/\epsilon) $ in the case of the unit disk it is sufficient.

2- Or, in the general case: we use a chart $ (f, \tilde \Omega) $ with $ f(0)=x_0 $ and we take $ \mu_{\epsilon}(x)= \eta_0 ( f( |x|/ \epsilon)) $ to have  connected  sets $ I_{\epsilon} $ and we take $ \eta_{\epsilon}(y)= \mu_{\epsilon}(f^{-1}(y))$. Because $ f, f^{-1} $ are Lipschitz, $ |f(x)-x_0| \leq k_ 2|x|\leq 1 $ for $ |x| \leq 1/k_2 $ and $ |f(x)-x_0| \geq k_ 1|x|\geq 2 $ for $ |x| \geq 2/k_1>1/k_2 $, the support  of $ \eta $ is in $ I_{(2/k_1)\epsilon} $.

$$ \begin{cases}
    
\eta_{\epsilon} \equiv 1,\,\,\,  {\rm on } \,\,\,  f(I_{(1/k_2)\epsilon}), \,\,\, 0 < \epsilon < \delta/2,\\

\eta_{\epsilon} \equiv 0,\,\,\, {\rm outside} \,\,\, f(I_{(2/k_1)\epsilon }), \\

0 \leq \eta_{\epsilon} \leq 1, \\

||\nabla \eta_{\epsilon}||_{L^{\infty}(I_{(2/k_1)\epsilon})} \leq \dfrac{C_0(\Omega, x_0)}{\epsilon}.

\end{cases} $$

3- Also, we can take: $ \mu_{\epsilon}(x)= \eta_0(|x|/\epsilon) $ and $ \eta_{\epsilon}(y)= \mu_{\epsilon}(f^{-1}(y)) $, we extend it by $ 0 $ outside $ f(B_1(0)) $.  We have $ f(B_1(0)) = D_1(x_0) $, $ f (B_{\epsilon}(0))= D_{\epsilon}(x_0) $ and $ f(B_{\epsilon}^+)= D_{\epsilon}^+(x_0) $ with $ f $ and $ f^{-1} $ smooth diffeomorphism.

$$ \begin{cases}
    
\eta_{\epsilon} \equiv 1,\,\,\,  {\rm on \, a \, the \, connected \, set } \,\,\,  J_{\epsilon} =f(I_{\epsilon}), \,\,\, 0 < \epsilon < \delta/2,\\

\eta_{\epsilon} \equiv 0,\,\,\, {\rm outside} \,\,\, J'_{\epsilon}=f(I_{2\epsilon }), \\

0 \leq \eta_{\epsilon} \leq 1, \\

||\nabla \eta_{\epsilon}||_{L^{\infty}(J'_{\epsilon})} \leq \dfrac{C_0(\Omega, x_0)}{\epsilon}.

\end{cases} $$

And, $ H_1(J'_{\epsilon}) \leq C_1 H_1(I_{2\epsilon}) = C_1 4\epsilon $, since $ f $ is Lipschitz. Here $ H_1 $ is the Hausdorff measure.

We solve the Dirichlet Problem:

\begin{displaymath}  \left \{ \begin {split} 
      -\Delta \bar \eta_{\epsilon}  & = -\Delta \eta_{\epsilon}              \,\, &&\text{in} \!\!&&\Omega \subset {\mathbb R}^2, \\
                  \bar \eta_{\epsilon} & = 0   \,\,             && \text{in} \!\!&&\partial \Omega.               
\end {split}\right.
\end{displaymath}

and finaly we set $ \tilde \eta_{\epsilon} =-\bar \eta_{\epsilon} + \eta_{\epsilon} $. Also, by the maximum principle and the elliptic estimates we have :

$$ ||\nabla \tilde \eta_{\epsilon}||_{L^{\infty}} \leq C(|| \eta_{\epsilon}||_{L^{\infty}} +||\nabla \eta_{\epsilon}||_{L^{\infty}} + ||\Delta \eta_{\epsilon}||_{L^{\infty}}) \leq \dfrac{C_1}{\epsilon^2}, $$

with $ C_1 $ depends on $ \Omega $.

We use the following estimate, see [8],

$$ ||\nabla v_i||_{L^q}\leq C_q,\,\, ||\nabla u_i||_q \leq C_q, \,\,\forall \,\, i\,\, {\rm and  }  \,\, 1< q < 2. $$

We deduce from the last estimate that, $ (v_i) $ converge weakly in $ W_0^{1, q}(\Omega) $, almost everywhere to a function $ v \geq 0 $ and $ \int_{\Omega} |x|^{2\beta} e^v < + \infty $ (by Fatou lemma). Also, $ V_i $ weakly converge to a nonnegative function $ V $ in $ L^{\infty} $. 

We deduce from the last estimate that, $ (u_i) $ converge weakly in $ W_0^{1, q}(\Omega) $, almost everywhere to a function $ u \geq 0 $ and $ \int_{\Omega} e^u < + \infty $ (by Fatou lemma). Also, $ W_i $ weakly converge to a nonnegative function $ W $ in $ L^{\infty} $.

The function $ u, v $ are in $ W_0^{1, q}(\Omega) $ solutions of :

\begin{displaymath} \left \{ \begin {split} 
      -\Delta u  & = |x|^{2 \beta} V e^{v} \in L^1(\Omega)    \,\, &&\text{in} \!\!&&\Omega \subset {\mathbb R}^2, \\
                  u  & = 0  \,\,                                     && \text{in} \!\!&&\partial \Omega.               
\end {split}\right.
\end{displaymath}
 
And,

\begin{displaymath} \left \{ \begin {split} 
      -\Delta v  & = W e^{u} \in L^1(\Omega)    \,\, &&\text{in} \!\!&&\Omega \subset {\mathbb R}^2, \\
                  v  & = 0  \,\,                                     && \text{in} \!\!&&\partial \Omega.               
\end {split}\right.
\end{displaymath}

 According to the corollary 1 of Brezis-Merle's result, see [6],   we have $ e^{k u }\in L^1(\Omega), k >1 $. By the elliptic estimates, we have $ v \in C^1(\bar \Omega) $.

According to the corollary 1 of Brezis-Merle's result, see [6],   we have $ e^{k v }\in L^1(\Omega), k >1 $. By the elliptic estimates, we have $ u \in C^1(\bar \Omega) $.

For two vectors $ f $ and $ g $ we denote by $  f \cdot g $ the inner product of $ f $ and $ g $. 

We can write:

\be -\Delta ((u_i-u) \tilde \eta_{\epsilon})= |x|^{2 \beta}(V_i e^{v_i} -Ve^v)\tilde \eta_{\epsilon} -2\nabla (u_i- u)\cdot \nabla \tilde \eta_{\epsilon}. \label{(1)}\ee

$$ -\Delta ((v_i-v) \tilde \eta_{\epsilon})= (W_i e^{u_i} -We^u)\tilde \eta_{\epsilon} -2\nabla (v_i- v)\cdot \nabla \tilde \eta_{\epsilon}. $$

We use the interior esimate of Brezis-Merle, see [6],

\bigskip

\underbar {\it Step 1:} Estimate of the integral of the first term of the right hand side of $ (\ref{(1)}) $.

\bigskip

We use the Green formula between $ \tilde \eta_{\epsilon} $ and $ u $, we obtain,

\be  \int_{\Omega} |x|^{2 \beta}Ve^v \tilde \eta_{\epsilon} dx =\int_{\partial \Omega} \partial_{\nu} u \eta_{\epsilon} \leq C'\epsilon ||\partial_{\nu}u||_{L^{\infty}}= C \epsilon \label{(2)}\ee

We have,

\begin{displaymath} \left \{ \begin {split} 
      -\Delta u_i  & = |x|^{2 \beta} V_i e^{v_i}                      \,\, &&\text{in} \!\!&&\Omega \subset {\mathbb R}^2, \\
                  u_i  & = 0  \,\,                                     && \text{in} \!\!&&\partial \Omega.               
\end {split}\right.
\end{displaymath}

We use the Green formula between $ u_i $ and $ \tilde \eta_{\epsilon} $ to have:

\be \int_{\Omega} |x|^{2 \beta} V_i e^{v_i} \tilde \eta_{\epsilon} dx = \int_{\partial \Omega} \partial_{\nu} u_i \eta_{\epsilon} d\sigma \to \mu_1(\eta_{\epsilon}) \leq \mu_1(J'_{\epsilon}) \leq 4  \pi - \epsilon_0, \,\,\, \epsilon_0 >0 \label{(3)}\ee

From $ (\ref{(2)}) $ and $ (\ref{(3)}) $ we have for all $ \epsilon >0 $ there is $ i_0 =i_0(\epsilon) $ such that, for $ i \geq i_0 $,

\be \int_{\Omega} ||x|^{2 \beta}(V_ie^{v_i}-Ve^v) \tilde \eta_{\epsilon}| dx \leq 4 \pi -\epsilon_0 + C \epsilon \label{(4)}\ee

\bigskip

\underbar {\it Step 2:} Estimate of integral of the second term of the right hand side of $ (\ref{(1)}) $.

\bigskip

Let $ \Sigma_{\epsilon} = \{ x \in \Omega, d(x, \partial \Omega) = \epsilon^3 \} $ and $ \Omega_{\epsilon^3} = \{ x \in \Omega, d(x, \partial \Omega) \geq \epsilon^3 \} $, $ \epsilon > 0 $. Then, for $ \epsilon $ small enough, $ \Sigma_{\epsilon} $ is hypersurface.

The measure of $ \Omega-\Omega_{\epsilon^3} $ is $ k_2\epsilon^3 \leq meas(\Omega-\Omega_{\epsilon^3})= \mu_L (\Omega-\Omega_{\epsilon^3}) \leq k_1 \epsilon^3 $.

{\bf Remark}: for the unit ball $ {\bar B(0,1)} $, our new manifold is $ {\bar B(0, 1-\epsilon^3)} $.

( Proof of this fact; let's consider $ d(x, \partial \Omega) = d(x, z_0), z_0 \in \partial \Omega $, this imply that $ (d(x,z_0))^2 \leq (d(x, z))^2 $ for all $ z \in \partial \Omega $ which it is equivalent to $ (z-z_0) \cdot (2x-z-z_0) \leq 0 $ for all $ z \in \partial \Omega $, let's consider a chart around $ z_0 $ and $ \gamma (t) $ a curve in $ \partial \Omega $, we have;

$ (\gamma (t)-\gamma(t_0) \cdot (2x-\gamma(t)-\gamma(t_0)) \leq 0 $ if we divide by $ (t-t_0) $ (with the sign and tend $ t $ to $ t_0 $), we have $ \gamma'(t_0) \cdot (x-\gamma(t_0)) = 0 $, this imply that $ x= z_0-s \nu_0 $ where $ \nu_0 $ is the outward normal of $ \partial \Omega $ at $ z_0 $))

With this fact, we can say that $ S= \{ x, d(x, \partial \Omega) \leq \epsilon \}= \{ x= z_0- s \nu_{z_0}, z_0 \in \partial \Omega, \,\, -\epsilon \leq s \leq \epsilon \} $. It  is sufficient to work on  $ \partial \Omega $. Let's consider a charts $ (z, D=B(z, 4 \epsilon_z), \gamma_z) $ with $ z \in \partial \Omega $ such that $ \cup_z B(z, \epsilon_z) $ is  cover of $ \partial \Omega $ .  One can extract a finite cover $ (B(z_k, \epsilon_k)), k =1, ..., m $, by the area formula the measure of $ S \cap B(z_k, \epsilon_k) $ is less than a $ k\epsilon $ (a $ \epsilon $-rectangle).  For the reverse inequality, it is sufficient to consider one chart around one point of the boundary.

We write,

\be \int_{\Omega} |\nabla ( u_i -u) \cdot \nabla \tilde \eta_{\epsilon} | dx =
\int_{\Omega_{\epsilon^3}} |\nabla (u_i -u) \cdot \nabla \tilde \eta_{\epsilon}| dx + \int_{\Omega - \Omega_{\epsilon^3}} |\nabla (u_i-u) \cdot \nabla \tilde \eta_{\epsilon}| dx.  \label{(5)}\ee

\bigskip

\underbar {\it Step 2.1:} Estimate of $ \int_{\Omega - \Omega_{\epsilon^3}} |\nabla (u_i-u) \cdot \nabla \tilde \eta_{\epsilon}| dx $.

\bigskip

First, we know from the elliptic estimates that  $ ||\nabla \tilde \eta_{\epsilon}||_{L^{\infty}} \leq C_1 /\epsilon^2 $, $ C_1 $ depends on $ \Omega $

We know that $ (|\nabla u_i|)_i $ is bounded in $ L^q, 1< q < 2 $, we can extract  from this sequence a subsequence which converge weakly to $ h \in L^q $. But, we know that we have locally the uniform convergence to $ |\nabla u| $ (by Brezis-Merle's theorem), then, $ h= |\nabla u| $ a.e. Let $ q' $ be the conjugate of $ q $.

We have, $  \forall f \in L^{q'}(\Omega)$

$$ \int_{\Omega} |\nabla u_i| f dx \to \int_{\Omega} |\nabla u| f dx $$

If we take $ f= 1_{\Omega-\Omega_{\epsilon^3}} $, we have:

$$ {\rm for } \,\, \epsilon >0 \,\, \exists \,\, i_1 = i_1(\epsilon) \in {\mathbb N}, \,\,\, i \geq  i_1,  \,\, \int_{\Omega-\Omega_{\epsilon^3}} |\nabla u_i| \leq \int_{\Omega-\Omega_{\epsilon^3}} |\nabla u| + \epsilon^3. $$

Then, for $ i \geq i_1(\epsilon) $,

$$ \int_{\Omega-\Omega_{\epsilon^3}} |\nabla u_i| \leq meas(\Omega-\Omega_{\epsilon^3}) ||\nabla u||_{L^{\infty}} + \epsilon^3 = \epsilon^3(k_1 ||\nabla u||_{L^{\infty}} + 1). $$

Thus, we obtain,

\be \int_{\Omega - \Omega_{\epsilon^3}} |\nabla (u_i-u) \cdot \nabla \tilde \eta_{\epsilon}| dx \leq  \epsilon C_1(2 k_1 ||\nabla u||_{L^{\infty}} + 1) \label{(6)}\ee

The constant $ C_1 $ does  not depend on $ \epsilon $ but on $ \Omega $.

\bigskip

\underbar {\it Step 2.2:} Estimate of $ \int_{\Omega_{\epsilon^3}} |\nabla (u_i-u) \cdot \nabla \tilde \eta_{\epsilon}| dx $.

\bigskip

We know that, $ \Omega_{\epsilon} \subset \subset \Omega $, and ( because of Brezis-Merle's interior estimates) $ u_i \to u $ in $ C^1(\Omega_{\epsilon^3}) $. We have,

$$ ||\nabla (u_i-u)||_{L^{\infty}(\Omega_{\epsilon^3})} \leq \epsilon^3,\, {\rm for } \,\, i \geq i_3 = i_3(\epsilon). $$

We write,
 
$$ \int_{\Omega_{\epsilon^3}} |\nabla (u_i-u) \cdot \nabla \tilde \eta_{\epsilon}| dx \leq ||\nabla (u_i-u)||_{L^{\infty}(\Omega_{\epsilon^3})} ||\nabla \tilde \eta_{ \epsilon}||_{L^{\infty}} \leq C_1 \epsilon \,\, {\rm for } \,\, i \geq i_3, $$

For $ \epsilon >0 $, we have for $ i \in {\mathbb N} $, $ i \geq \max \{i_1, i_2, i_3 \} $,

\be \int_{\Omega} |\nabla (u_i-u) \cdot \nabla \tilde \eta_{\epsilon}| dx \leq \epsilon C_1(2 k_1 ||\nabla u||_{L^{\infty}} + 2) \label{(7)}\ee

From $ (\ref{(4)}) $ and $ (\ref{(7)}) $, we have, for $ \epsilon >0 $, there is $ i_3= i_3(\epsilon) \in {\mathbb N}, i_3 = \max \{ i_0, i_1, i_2 \} $ such that,

\be \int_{\Omega} |-\Delta [(u_i-u)\tilde \eta_{\epsilon}]|dx \leq 4 \pi-\epsilon_0+  \epsilon 2 C_1(2 k_1 ||\nabla u||_{L^{\infty}} + 2 + C) \label{(8)}\ee

We choose $ \epsilon >0 $ small enough to have a good estimate of  $ (\ref{(1)}) $.

Indeed, we have:

\begin{displaymath} \left \{ \begin {split} 
      -\Delta [(u_i-u) \tilde \eta_{\epsilon}]   & = g_{i,\epsilon}                   \,\, &&\text{in} \!\!&&\Omega \subset {\mathbb R}^2, \\
                 (u_i-u) \tilde \eta_{\epsilon}    & = 0  \,\,                                     && \text{in} \!\!&&\partial \Omega.               
\end {split}\right.
\end{displaymath}

with $ ||g_{i, \epsilon} ||_{L^1(\Omega)} \leq 4 \pi -\dfrac{\epsilon_0}{2}. $

We can use Theorem 1 of [6] to conclude that there are $ q \geq \tilde q >1 $ such that:

$$ \int_{V_{\epsilon}(x_0)} e^{\tilde q |u_i-u|} dx \leq \int_{\Omega} e^{q|u_i -u| \tilde \eta_{\epsilon}} dx \leq C(\epsilon,\Omega). $$
 
where, $ V_{\epsilon}(x_0) $ is a neighberhood of $ x_0 $ in $ \bar \Omega $. Here we have used that in a neighborhood of $ x_0 $  by the elliptic estimates, 
$ 1- C \epsilon \leq \tilde \eta_{\epsilon} \leq 1 $.

Thus, for each $ x_0 \in \partial \Omega - \{ \bar x_1,\ldots, \bar x_m \} $ there is $ \epsilon_{x_0} >0, q_{x_0} > 1 $ such that:

\be \int_{B(x_0, \epsilon_{x_0})} e^{q_{x_0} u_i} dx \leq C, \,\,\, \forall \,\,\, i. \label{(9)}\ee

Now, we consider a cutoff function $ \eta \in C^{\infty}({\mathbb R}^2) $ such that

$$ \eta \equiv 1 \,\,\, {\rm on } \,\,\, B(x_0, \epsilon_{x_0}/2) \,\,\, {\rm and } \,\,\, \eta \equiv 0 \,\,\, {\rm on } \,\,\, {\mathbb R}^2 -B(x_0, 2\epsilon_{x_0}/3). $$

We write

$$ -\Delta (v_i \eta) = W_i e^{u_i} \eta - 2 \nabla v_i \cdot \nabla \eta  - v_i \Delta \eta. $$

Because, by Poincar\'e and Gagliardo-Nirenberg-Sobolev inequalities:

$$ ||v_i||_{q^*} \leq c_q ||\nabla v_i||_q \leq  C_q, \,\, 1\leq q <2, $$

with, $ q^*=2q/(2-q) >2 >1 $.

\smallskip

By the elliptic estimates, $ (v_i)_i $ is uniformly bounded in $ L^{\infty}(V_{\epsilon}(x_0)) $. Finaly, we have, for some $ \epsilon > 0 $ small enough,

$$ || v_i||_{C^{0,\theta}[B(x_0, \epsilon)]} \leq c_3 \,\,\, \forall \,\,\, i. $$

Now, we consider a cutoff function $ \eta \in C^{\infty}({\mathbb R}^2) $ such that

$$ \eta \equiv 1 \,\,\, {\rm on } \,\,\, B(x_0, \epsilon_{x_0}/2) \,\,\, {\rm and } \,\,\, \eta \equiv 0 \,\,\, {\rm on } \,\,\, {\mathbb R}^2 -B(x_0, 2\epsilon_{x_0}/3). $$

We write

$$ -\Delta (u_i \eta) = |x|^{2\beta} V_i e^{v_i} \eta - 2 \nabla u_i \cdot \nabla \eta  - u_i \Delta \eta. $$

By the elliptic estimates, $ (u_i)_i $ is uniformly bounded in $ L^{\infty}(V_{\epsilon}(x_0)) $ and also in $ C^{0,\theta} $ norm.

If we repeat this procedure another time, we have a boundedness of $ (u_i)_i $ and $ (v_i)_i $ in the $ C^{1,\theta} $ norm, because they are bounded in $ W^{2,q}\subset W^{1,q^*} $ norms with $ 2q/(2-q)=q^* >2 $.

We have the same computations and conclusion if we consider a regular point $ x_0 $ for the measure $ \mu_2 $.

We have proved that, there is a finite number of points $ \bar x_1, \ldots, \bar x_m $ such that the squence $ (u_i)_i  $ and $ (v_i)_i $ are locally uniformly bounded (in $ C^{1,\theta}, \theta >0 $) in $ \bar \Omega - \{ \bar x_1, \ldots , \bar x_m \} $.

\bigskip

\underbar {\it Proof of theorem 1.2:} 

Without loss of generality, we can assume that $ 0=x_1 $ is a blow-up point. Since the boundary is an analytic curve $ \gamma_1(t) $, there is a neighborhood of  $ 0 $ such that the curve $ \gamma_1 $ can be extend to a holomorphic map such that $ \gamma_1'(0) \not = 0 $ (series) and by the inverse mapping one can assume that this map is univalent around $ 0 $. In the case when the boundary is a simple Jordan curve the domain is simply connected. In the case that the domains has a finite number of holes it is conformally equivalent to a disk with a finite number of disks removed. Here we consider a general domain. Without loss of generality one can assume that $ \gamma_1 (B_1^+) \subset \Omega $ and also $ \gamma_1 (B_1^-) \subset (\bar \Omega)^c $ and $ \gamma_1 (-1,1) \subset \partial \Omega $ and $ \gamma_1 $ is univalent. This means that $ (B_1, \gamma_1) $ is a local chart around $ 0 $ for $ \Omega $ and $ \gamma_1 $ univalent. (This fact holds if we assume that we have an analytic domain, (below a graph of an analytic function), we have necessary the condition $ \partial \bar \Omega = \partial \Omega $ and the graph is analytic, in this case $ \gamma_1 (t)= (t, \phi(t)) $ with $ \phi $ real analytic and an example of this fact is the unit disk  around the point $ (0,1) $ for example).

By this conformal transformation, we can assume that $ \Omega =B_1^+ $, the half ball, and $ \partial^+ B_1^+ $ is the exterior part, a part which not contain $ 0 $ and on which  $ u_i $ converge in the $ C^1 $ norm to $ u $. Let us consider $ B_{\epsilon}^+ $, the half ball with radius $ \epsilon >0 $. Also, one can consider a $ C^1 $ domain (a rectangle between two half disks) and by charts its image is a $ C^1 $ domain)
We know that:

$$ u_i \in C^{2, \epsilon}(\bar \Omega). $$ 

Thus we can use integrations by parts (Stokes formula). The Pohozaev identity applied around the blow-up $ 0 $:

\be  \int_{B_{\epsilon}^+} \Delta u_i < x |\nabla v_i > dx = -  \int_{B_{\epsilon}^+} \Delta v_i < x |\nabla u_i > dx + \int_{\partial^+ B_{\epsilon}^+}  g(\nabla u_i, \nabla v_i)d\sigma, \label{(10)}\ee

Thus,

\be  \int_{B_{\epsilon}^+} |x|^{2 \beta} V_i e^{v_i}< x |\nabla v_i > dx = -  \int_{B_{\epsilon}^+} W_ie^{u_i} < x |\nabla u_i > dx - \int_{\partial^+ B_{\epsilon}^+}  g(\nabla u_i, \nabla v_i)d\sigma, \label{(11)}\ee

After integration by parts, we obtain:

$$  \int_{B_{\epsilon}^+} 2V_i(1+\beta) |x|^{2 \beta} e^{v_i} dx +  \int_{B_{\epsilon}^+} < x |\nabla V_i > |x|^{2 \beta} e^{v_i} dx+  \int_{\partial B_{\epsilon}^+} < \nu |x >|x|^{2 \beta} V_i e^{v_i} d\sigma+ $$

$$ + \int_{B_{\epsilon}^+} W_i e^{u_i} dx +  \int_{B_{\epsilon}^+} < x |\nabla W_i > e^{u_i} dx+  \int_{\partial B_{\epsilon}^+} < \nu | x > W_i e^{u_i}  d\sigma = $$

$$ = - \int_{\partial^+ B_{\epsilon}^+} g(\nabla u_i, \nabla v_i)d\sigma, \label{(10)} $$ 

Also, for $ u $ and $ v $, we have:

$$  \int_{B_{\epsilon}^+} 2V(1+ \beta) |x|^{2 \beta} e^{v} dx +  \int_{B_{\epsilon}^+} < x |\nabla V > |x|^{2 \beta} e^{v} dx+  \int_{\partial B_{\epsilon}^+} < \nu | x >|x|^{2 \beta} V e^{v}  d\sigma+ $$

$$ + \int_{B_{\epsilon}^+} We^{u} dx +  \int_{B_{\epsilon}^+} < x |\nabla W> e^{u} dx+  \int_{\partial B_{\epsilon}^+} < \nu | x > W e^{u} d\sigma = $$

$$ = - \int_{\partial^+ B_{\epsilon}^+} g(\nabla u, \nabla v)d\sigma, $$ 

If, we take the difference, we obtain:

$$ \int_{\gamma_1 (B_{\epsilon}^+)} |x|^{2 \beta} V_i e^{v_i} dx + \int_{\gamma_1 (B_{\epsilon}^+)} W_i e^{u_i} dx  = o(\epsilon)+o(1) $$

But,

 $$ \int_{\gamma_1 (B_{\epsilon}^+)} |x|^{2 \beta} V_i e^{v_i} dx + \int_{\gamma_1 (B_{\epsilon}^+)} W_i e^{u_i} dx  = \int_{\partial \gamma_1 (B_{\epsilon}^+)} \partial_{\nu} u_i d\sigma + \int_{\partial \gamma_1 (B_{\epsilon}^+)} \partial_{\nu} v_i d\sigma $$

and,

\be \int_{\partial \gamma_1 (B_{\epsilon}^+)} \partial_{\nu} u_i d\sigma + \int_{\partial \gamma_1 (B_{\epsilon}^+)} \partial_{\nu} v_i d\sigma 
\to \alpha_1+\beta_1 >0 \label{(12)}\ee

a contradiction.

\end{document}